\documentclass[a4paper]{amsart}

\usepackage{amsthm}
\usepackage{amscd}
\usepackage{graphics}
\usepackage{amssymb,amsmath}  
\usepackage{latexsym} 
\usepackage[english]{babel}
\usepackage{tabularx}
\usepackage{longtable}

\usepackage{a4wide}

\newtheorem{lem}{Lemma}[section]
    \newtheorem{prop}[lem]{Proposition}
    \newtheorem{thm}[lem]{Theorem}

   \theoremstyle{definition}
    \newtheorem{dfn}[lem]{Definition}

\theoremstyle{remark}

\numberwithin{equation}{section}

\DeclareMathOperator{\rank}{rank}

\def\bb{\mathbf}
\def\phi{\varphi}
\newcommand{\Pp}[1]{\bb{P}^{#1}}
\newcommand{\GL}{\mathrm{GL}}
\newcommand{\PGL}{\mathrm{PGL}}

\def\Z{\bb{Z}}
\def\Q{\bb{Q}}
\def\C{\bb{C}}
\def\co{\colon\thinspace}
\def\pu{\bullet}
\def\I{\mathcal I}

\newcommand{\Mm}[2]{\mathcal M_{{#1},{#2}}}

\newcommand{\Hhm}[2]{\mathcal H_{{#1},{#2}}}
\newcommand{\s}{\mathfrak S}

\DeclareMathOperator{\Aut}{Aut}
\DeclareMathOperator{\Phg}{\wp}
\newcommand{\Phgs}[1][{\s_2}]{\Phg^{#1}}
\newcommand{\Ll}{\mathbf L}

\newcommand{\Schur}[1]{\bb S_{#1}}
\def\Inv{\Schur 2}
\def\Ant{\Schur{1,1}}
\newcommand{\schur}[1]{s_{#1}}
\def\inv{\schur 2}
\def\ant{\schur{1,1}}

\newcommand{\coh}[3][\Q]{H^{#2}({#3};{#1})}

\DeclareMathOperator{\Fil}{Fil}
\def\ba{\big|}
\def\X{\mathcal{X}}
\DeclareMathOperator{\Conf}{Conf}
\newcommand{\B}[2]{B({#2},{#1})}
\newcommand{\tB}[2]{\tilde{B}({#2},{#1})}
\newcommand{\F}[2]{F({#2},{#1})}
\newcommand{\tF}[2]{\tilde{F}({#2},{#1})}

\def\op{\mathring}
\def\duale{\check{\ }}

\begin{document}

\title{Rational cohomology of $\Mm32$}
\date{February 7th, 2007.}

\author{Orsola Tommasi}
\email{tommasi@mathematik.uni-mainz.de}
\address{Institut f\"ur Mathematik,
FB 08 --- Physik, Mathematik und Informatik,
Johannes Gutenberg-U\-ni\-ver\-si\-t\"at,
D-55099 Mainz,
Germany}
\curraddr{Institut Mittag-Leffler,
Aurav\"agen 17,
SE-182 60 Djursholm,
Sweden}

\subjclass[2000]{14H10 (primary), 55R80, 14F99 (secondary).}
\keywords{Moduli space, $2$-pointed curves of genus three, rational cohomology, discriminants.}
\thanks{Support from the Mittag-Leffler Institute during the last part of the preparation of this paper is gratefully acknowledged.}

\begin{abstract}
We compute the rational cohomology of the moduli space of non-singular complex plane quartic curves with two marked points. This allows to calculate the rational cohomology of the moduli space of non-singular complex projective curves of genus~3 with two marked points.
\end{abstract}

\maketitle

\section{Introduction}\label{introM3n}

Let us denote by $\Mm32$ the moduli space of non-singular complex projective curves with two marked points, and by $\mathcal Q_2$ the moduli space of plane quartic curves with two marked points. In this paper, we prove
 
\begin{thm}\label{result}
The rational cohomology groups of $\mathcal Q_2$ and $\Mm32$, with their mixed Hodge structures and their structures as $\s_2$-representations, are as follows.
\begin{enumerate}
\item $\coh k {\mathcal Q_2}=\left\{
\begin{array}{ll}
\Inv\otimes\Q &k=0,\\
(\Inv+\Ant)\otimes\Q(-1) &k=2,\\
\Inv\otimes\Q(-3) &k=5,\\
\Inv\otimes\Q(-6) &k=6,\\
(\Inv+\Ant)\otimes\Q(-7)+\Ant\otimes\Q(-8) &k=8,\\
0 &\text{otherwise}.\\
\end{array}
\right.$
\item\label{resultii}  $\coh k {\Mm32}=\left\{
\begin{array}{ll}
\Inv\otimes\Q &k=0,\\
(\bigoplus^2\Inv+\Ant)\otimes\Q(-1) &k=2,\\
(\Inv+\Ant)\otimes\Q(-2) &k=4,\\
\Inv\otimes\Q(-3) &k=5,\\
\Inv\otimes\Q(-6) &k=6,\\
(\Inv+\Ant)\otimes\Q(-7) &k=8,\\
0 &\text{otherwise}.\\
\end{array}
\right.$
\end{enumerate}
\end{thm}

Note that the Euler characteristic of $\Mm32$ in the Grothendieck group of mixed Hodge structures was computed in \cite{BT}, exploiting Jonas Bergstr\"om's count of the number of points of $\Mm32$ defined over finite fields (\cite{Jonas2}). Another computation of this Euler characteristic can be found in \cite[Chapter~III]{OTtesi}. All these results agree with the topological Euler characteristic of $\Mm32$ as calculated in~\cite{Eulerchar}.

It is well known that the canonical model of a non-hyperelliptic curve of genus $3$ is a smooth quartic in the projective plane. Hence $\mathcal Q_2$ is the complement of the hyperelliptic locus $\Hhm32$ inside $\Mm 32$.
Let us start by considering the moduli space $\mathcal Q$ of smooth quartic curves in the projective plane.
Quartic curves are defined by the vanishing of polynomials of degree four in three indeterminates, i.e., by elements of $S^2_4:=\C[x_0,x_1,x_2]_4$. Clearly, not every element of $S^2_4$ defines a non-singular curve, but we have to exclude the locus $\Sigma^2_4\subset S^2_4$ of singular polynomials. The action of $G=\GL(3)$ on the coordinates $x_0,x_1,x_2$ induces an action on $S^2_4\setminus\Sigma^2_4$, and $\mathcal Q$ is the the geometric quotient of $S^2_4\setminus\Sigma^2_4$ by the action of $G$.

The rational cohomology of $S^2_4\setminus\Sigma^2_4$ was computed by Vassiliev in \cite{Vart}. Comparing this result with the rational cohomology of the moduli space $\mathcal Q$, as computed by Looijenga in \cite{Looij}, one observes that the cohomology of the space of non-singular polynomials in $S^2_4$ is isomorphic (as graded vector space) to the tensor product of the cohomology of the moduli space  $\mathcal Q$ and that of $G=\GL(3)$. Indeed, Peters and Steenbrink \cite{PS} proved that this is always the case when comparing the rational cohomology of the space of non-singular homogeneous polynomials with the cohomology of the corresponding moduli space of smooth hypersurfaces. 

As explained in \cite[\S~5]{BT}, Peters--Steenbrink's result can be adapted to moduli spaces of smooth hypersurfaces with $m$ marked points, when $m$ is small enough. This requires to replace the space $S^n_d$ of homogeneous polynomials of degree $d$ in $x_0,\dots,x_n$ with a certain incidence correspondence. In our case ($n=2, d=4, m=2$) we set 
$$\I_2:=\{(\alpha,\beta,f)\in\F 2{\Pp2} \times(S^2_4\setminus\Sigma^2_4): f(\alpha)=f(\beta)=0\},$$
where $\F2{\Pp2}$ denotes the complement of the diagonal in $\Pp2\times\Pp2$.
The action of $G=\GL(3)$ on $\Pp2$ and $S^2_4$ can be extended to $\I_2$, and the geometric quotient $\I_2/G$ is isomorphic to $\mathcal Q_2$. Then the following isomorphism of graded vector spaces with mixed Hodge structures holds:
\begin{equation}\label{icLH}\coh\pu{\I_2}\cong \coh\pu{\mathcal Q_2}\otimes \coh\pu{\GL(3)}.
\end{equation}
This follows from \cite{PS}, in view of \cite[Theorem~5.2]{BT}.
As a consequence, we have that determining the rational cohomology of $\I_2$ immediately yields the rational cohomology of $\mathcal Q_2$. Note that the isomorphism~\eqref{icLH} is compatible with the action of the symmetric group $\s_2$ on the cohomology groups of $\mathcal Q_2$ and $\I_2$ induced by the involution interchanging the two marked points.

We compute $\coh\pu{\I_2}$ by studying the natural projection $\pi_2\co \I_2\rightarrow \F 2{\Pp2}$. 
The map $\pi_2$ is a locally trivial fibration, whose fibre is the complement of $\Sigma^2_4$ in a linear subspace of $S^2_4$. Therefore, we can compute the cohomology of this fibre with Vassiliev--Gorinov's method for the cohomology of complements of discriminants. The study of the Leray spectral sequence associated to the fibration $\pi_2$ allows to determine the cohomology of $\I_2$. In this last step, we will use very often the relation~\eqref{icLH}.

The plan of the paper is as follows. In \S\S~\ref{VGI} and~\ref{QKQK} we compute the rational cohomology of $\mathcal Q_2$ and prove Theorem~\ref{result} by the methods explained above. We conclude the paper with a concise review of Vassiliev--Gorinov's method.

\subsection*{Acknowledgements}
The author is grateful to Joseph Steenbrink, Alexei Gorinov, Carel Faber, Torsten Ekedahl  and Jonas Berg\-str\"om for useful discussions during different phases of the preparation of this paper. The author would also like to thank the referee for his careful work and especially for pointing out a mistake in an earlier version of part \eqref{resultii} of Theorem~\ref{result}.

\subsection*{Notation}

{\small
\noindent\begin{longtable}{ll}
$\s_n$&the symmetric group in $n$ letters.\\
$S^n_d$& vector space of homogeneous polynomials of degree $d$ in $n+1$ indeterminates $x_0,\dots,x_n$.\\
$\Sigma^n_d$&locus of singular polynomials in $S^n_d$.\\
$K_0(\mathsf{HS_\Q})$& Grothendieck group of rational (mixed) Hodge structures over $\Q$.\\
$K_0(\mathsf{HS_\Q^{\s_n}})$& Grothendieck group of rational (mixed) Hodge structures endowed with an $\s_n$-action.\\
$\Q(m)$& Tate Hodge structure of weight $-2m$.\\
$\Ll$& class of $\Q(-1)$ in $K_0(\mathsf{HS_\Q})$.\\
$\Schur\lambda$& $\Q$-representation of $\s_n$ indexed by the partition $\lambda\vdash n$.\\
$\schur\lambda$& Schur polynomial indexed by the partition $\lambda\vdash n$.\\
$\Delta_j$ & $j$-dimensional closed simplex.\\
$\op\Delta_j$ & interior of the $j$-dimensional closed simplex.\\
$\F kZ$& space of ordered configurations of $k$ distinct points on the variety $Z$ (see Def.~\ref{tuples}).\\
$\B kZ$& space of unordered configurations of $k$ distinct points on the variety $Z$ (see Def.~\ref{tuples}).
\end{longtable}}

Throughout this paper we will make an extensive use of Borel--Moore homology, i.e., homology with locally finite support. A reference for its definition and the properties we use is for instance \cite[Chapter 19]{Fulton}.

To write the results on cohomology and Borel--Moore homology groups in a compact way, we will express them by means of polynomials, in the following way. Let $T_\pu$ denote a graded $\Q$-vector space with mixed Hodge structures. For every $i\in\Z$, we can consider the class $[T_i]$ in the Grothendieck group of rational Hodge structures. We define the Hodge--Grothendieck  polynomial (for short, HG polynomial) of $T_\pu$ to be the polynomial
$$\Phg(T_\pu):=\sum_{i\in\Z} [T_i] t^i\in K_0(\mathsf{HS_\Q})[t].$$

If moreover a symmetric group $\s_n$ acts on $T_\pu$ respecting the grading and the mixed Hodge structures on $T_\pu$, we define the \emph{$\s_n$-equivariant} Hodge Grothendieck polynomial (for short, $\s_n$-HG polynomial) $\Phgs[\s_n](T_\pu)$ by replacing $ K_0(\mathsf{HS_\Q})$ by $ K_0(\mathsf{HS_\Q^{\s_n}})$ in the definition of the HG polynomial. 

\section{Rational cohomology of $\mathcal Q_2$}\label{VGI}

Consider the space $S^2_4$ of homogeneous quartic polynomials in $x_0,x_1,x_2$,  and denote by $\Sigma:=\Sigma^2_4$ the discriminant, i.e., the locus of singular quartic polynomials. For every $p\in\Pp2$, denote by $V_p$ the linear subspace of $S^2_4$ of polynomials vanishing at $p$. 
The aim of this section is to calculate the rational cohomology of the incidence correspondence 
$\I_2=\{(\alpha,\beta,f)\in\F2{\Pp2}\times(S^2_4\setminus\Sigma): f(\alpha)=f(\beta)=0\}.$
Note that knowing the cohomology of $\I_2$ is equivalent to knowing the cohomology of its projectivization
$$\mathcal P_2=\{(\alpha,\beta,[f])\in\F2{\Pp2}\times\Pp{}(S^2_4\setminus\Sigma): f(\alpha)=f(\beta)=0\},$$
as the rational cohomology of $\I_2$ is isomorphic to the tensor product of the cohomology of $\mathcal P_2$ and $\coh\pu{\C^*}$.

We will start by applying Vassiliev--Gorinov's method (see \S~\ref{VGmethod}) to the calculation of the cohomology of $(V_p\cap V_q)\setminus\Sigma$, where $p$ and $q$ are two fixed distinct points in $\Pp2$.  Next, we will consider the Leray spectral sequence for the natural projection $\pi_2\co  {\I_2}\rightarrow \F2{\Pp2}$. Note that the map $\pi_2$ is a locally trivial fibration with fibre isomorphic to $V_p\cap V_q \setminus\Sigma$.

By Alexander's duality between reduced cohomology and Borel--Moore homology, we have
\begin{equation}\label{alex}\tilde H^\pu((V_p\cap V_q)\setminus\Sigma;\Q)\cong\bar H_{25-\pu}(V_p\cap V_q\cap\Sigma;\Q)(-13).\end{equation}

To apply Vassiliev--Gorinov's method to $V_p\cap V_q\cap\Sigma$, we need an ordered list of all possible singular sets of the elements in $V_p\cap V_q\cap\Sigma$. 
We can easily obtain such a list by an adaptation of the list of possible singular configurations of quartic curves (like the one in \cite[Proposition~6]{Vart}). For every configuration in the list, one has to distinguish further whether the singular points are or are not in general position with respect to $p$ and $q$  (for instance, if $p$ or $q$ are or are not contained in the singular configuration). This procedure yields a complete list of singular sets of elements of $V_p\cap V_q$; let us denote by $R$ the number of types of configurations in the list. 
As recalled in \S~\ref{VGmethod}, Vassiliev--Gorinov's method gives a recipe to construct spaces $\ba\X\ba$, $\ba\Lambda\ba$ and a map
$$|\epsilon|\co  \ba \X\ba \longrightarrow V_p\cap V_q\cap\Sigma$$
inducing an isomorphism on Borel--Moore homology.
The Borel--Moore homology of $\ba \X\ba$ (respectively, $\ba \Lambda\ba$) can be computed by considering the stratification $\{F_j\}_{j=1,\dots,R}$ (resp., $\{\Phi_j\}_{j=1,\dots,R}$). The properties of $F_j$ and $\Phi_j$ are explained in Proposition~\ref{ucci}. Recall in particular that $F_j$ is the total space of a vector bundle over $\Phi_j$, and that for finite configurations the Borel--Moore homology of $\Phi_j$ coincides (after a shift in the indices) with the Borel--Moore homology of the space of configurations of type $j$ with coefficients in a rank $1$ local system changing its orientation every time two points in a configuration are interchanged.

In our case, for most indices $j\in\{1,\dots,R\}$ the space of singular configurations of type $j$ has trivial Borel--Moore homology in the appropriate system of coefficients. Hence, the strata $F_j$ have trivial Borel--Moore homology. 
In view of Lemma~\ref{lem1}, this is the case for configurations with too many points lying on the same rational curve. Furthermore, the same occurs for configurations containing rational curves as components (see \cite[Lemma~2.17]{OTM4} and following remarks). 

In Table~\ref{talist} we list all remaining configurations, i.e., all singular configurations indexing strata that give a non-trivial contribution to the Borel--Moore homology of $\Sigma\cap V_p\cap V_q$. In the same table, we also give a description of the strata of $\ba\Lambda\ba$ and $\ba\X\ba$ corresponding to each configuration. 
From the descriptions, it is straightforward to compute the Borel--Moore homology of the strata $\Phi_j$ and $F_j$ for $1\leq j\leq7$. 
The most difficult strata (corresponding to configurations of type 8, 9 and 10) are studied separately in \S~\ref{QKQK}. The results there, together with the description of the strata given in Table~\ref{talist}, allow to compute the $E^1$ terms of the spectral sequences in Borel--Moore homology converging to $\ba\Lambda\ba$ and $\ba\X\ba$, induced by the filtrations associated, respectively, with $\{\Phi_j\}$ and $\{F_j\}$.

\begin{table}\caption{\label{talist} Singular configurations and their contribution}
\begin{tabular}{l@{}p{14cm}}
1.&The point $p$ or the point $q$. \\
&Stratum:  $F_1$ is a $\C^{11}$-bundle over $\Phi_1=\{p,q\}$.\\
2. &Any point different from $p,q$. \\
&Stratum: $F_2$ is a $\C^{10}$-bundle over $\Phi_2\cong\Pp2\setminus\{p,q\}$.\\
3.&The pair $\{p,q\}$. \\
&Stratum:  $F_3$ is a $\C^9$-bundle over $\Phi_3\cong\op\Delta_1$.\\
4a.\;&Pairs of points on the line $pq$, different from $\{p,q\}$.\\
&Stratum:  $F_{4a}$ is a $\C^8$-bundle over $\Phi_{4a}$, which is a non-orientable $\op\Delta_1$-bundle over a space which can be decomposed as the disjoint union of $\C^*$ and $\B2\C$.\\
4b.\;& Pairs of points $\{a,b\}$ with $a\in\{p,q\}$, $b\notin pq$. \\
&Stratum:  $F_{4b}$ is a $\C^8$-bundle over $\Phi_{4b}$, which is  a non-orientable $\op\Delta_1$-bundle over the disjoint union of two copies of $\C^2$.\\
5.&Pairs of points $\{a,b\}$ with $a\in(pq\setminus\{p,q\})$, $b\notin pq$. \\
&Stratum:  $F_5$ is a $\C^7$-bundle over $\Phi_5$, which is  a non-orientable $\op\Delta_1$-bundle over $\C^*\times\C^2$.\\
6.&Triplets consisting of $p,q$ and another point outside $pq$. \\
&Stratum:  $F_6$ is a $\C^6$-bundle over $\Phi_6$, which is  a $\op\Delta_2$-bundle over $\C^2$.\\
7.&Triplets with two points on $pq$ (not both in $\{p,q\}$) and another point outside $pq$. \\
&Stratum:  $F_7$ is a $\C^5$-bundle over $\Phi_7$, which is  a non-orientable $\op\Delta_2$-bundle over a space that can be decomposed as the disjoint union of $\C^*\times\C^2$ and $\B2\C\times\C^2$.\\
{8}.& Five points $a,b,c,d,e\in\Pp2$, such that $a,b,d,e,p,q$ lie on a conic different from $ab\cup de$, $\{c\}=ab\cap de\not\subset\{p,q\}$ and $\{p,q\}\not\subset\{a,b,d,e\}$.\\
&Stratum:  $F_{8}$ is a $\C$-bundle over $\Phi_{8}$, which is  a $\op\Delta_4$-bundle over the configuration space $X_{8}$ of \S~\ref{QKQK}.
\\
{9}.& Six points that are the pairwise intersection of four lines $\ell_i$ ($1\leq i\leq4$) in general position, such that $\{p,q\}\subset\bigcup_i\ell_i.$ \\
&Stratum:  
$F_{9}$ is a $\C$-bundles over $\Phi_{9}$, which is a $\op\Delta_5$-bundle over the configuration spaces $X_{9}$ studied in \S~\ref{QKQK}. The simplices bundle does not change its orientation when two lines $\ell_i,\ell_j$ are interchanged.\\
{10}.& The entire $\Pp2$. \\
&Stratum:  $F_{10}$ is an open cone over the space $\ba\Lambda\ba$, which is the union of all strata $\Phi_j$ with $j\leq 9$.
\end{tabular}
\end{table}

The columns of the spectral sequences converging to the Borel--Moore homology of $\ba\X\ba$ and $\ba\Lambda\ba$ can be divided into two blocks: one with the first seven columns, the other with columns 8, 9 and 10. Looking at Hodge weights, one can easily prove that all differentials in the spectral sequence between columns in the block 1--7 and in the block 8--10 are trivial. Furthermore, this behaviour carries on when one investigates the Leray spectral sequence associated to the fibration $\pi_2$. Therefore, we will consider the two blocks separately. The contribution of columns 1--7 is computed below. The contribution of columns 8--10 to the rational cohomology of $\I_2$ is computed in \S~\ref{QKQK}.

In the spectral sequence converging to the Borel--Moore homology of $\ba\Lambda\ba$ all terms in the first seven columns are killed by differentials, with the exception of an $\s_2$-invariant 1-dimensional homology group in degree $0$. This follows from dimensional reasons: If these classes were not killed, they would give rise to cohomology classes of degree $\geq 14$ in the cohomology of $(V_p\cap V_q)\setminus\Sigma$, and this is impossible because the latter is affine of dimension 13. As a consequence, the strata 1--7 do not contribute to the Borel--Moore homology of the open cone $F_{10}$.

The first seven columns of the spectral sequence converging to the Borel--Moore homology of $\Fil_j\ba\X\ba$  are given in Table~\ref{t2}. Note that the description of the strata of the domain $\ba\X\ba$ of the geometric realization given in Table~\ref{talist} allows us to study the behaviour of each Borel--Moore homology class with respect to the interchange of the points $p,q$. Table~\ref{t2} includes also the information on the $\s_2$-action generated by this involution.

\begin{table}
\caption{\label{t2}
First seven columns of the spectral sequence converging to the Borel--Moore homology of $V_p\cap V_q\cap\Sigma$}
\scalebox{0.89}{
$\begin{array}{r|c@{\,}c@{\,}c@{\,}c@{\,}c@{\,}c@{\,}c}
22&0&\Inv\otimes\Q(12)&0&0&0&0&0\\
21&(\Inv+\Ant)\otimes\Q(11)&0&0&0&0&0&0\\
20&0&\Inv\otimes\Q(11)&0&0&0&0&0\\
19&0&\Ant\otimes\Q(10)&0&0&0&0&0\\
18&0&0&0&0&0&0&0\\
17&0&0&0&(\Inv+\Ant)\otimes\Q(10)&0&0&0\\
16&0&0&\Ant\otimes\Q(9)&0&\Inv\otimes\Q(10)&0&0\\
15&0&0&0&\Inv\otimes\Q(9)&\Ant\otimes\Q(9)&0&0\\
14&0&0&0&\Ant\otimes\Q(8)&0&0&0\\
13&0&0&0&0&0&0&0\\
12&0&0&0&0&0&\Ant\otimes\Q(8)&0\\
11&0&0&0&0&0&0&\Inv\otimes\Q(8)\\
10&0&0&0&0&0&0&\Ant\otimes\Q(7)\\\hline
 &1&2&3&4&5&6&7
\end{array}
$
}
\end{table}

In the spectral sequence in Table~\ref{t2}, the only possibly non-trivial differential is $d_2\co E^2_{5,15}\rightarrow E^2_{3,16}$. This is certainly zero, because otherwise we would get a contradiction with the isomorphism 
$$\coh\pu{(V_p\cap V_q)\setminus\Sigma}\cong\coh\pu{\C^*}\otimes\coh\pu{\Pp{}((V_p\cap V_q)\setminus\Sigma)}.$$ 
In view of Alexander's duality \eqref{alex}, we have that the part of the rational cohomology of $V_p\cap V_q\setminus\Sigma$ that comes from the first seven columns of Vassiliev--Gorinov's spectral sequence has $\s_2$-HG polynomial
\begin{equation}\label{1to7}
(1+\Ll t)\big(\inv+\Ll^2t^3(2\inv+\ant)+\Ll^3t^4\ant+\Ll^4t^6(\inv+\ant)+\Ll^5t^7\ant\big),
\end{equation}
where the $\s_2$-action is generated by the involution interchanging $p$ and $q$. Note that the second factor of \eqref{1to7} is the $\s_2$-HG polynomial of the cohomology of the projectivization of $V_p\cap V_q\setminus\Sigma$.

Next, we study the contribution of this part of the cohomology of $V_p\cap V_q\setminus\Sigma$ to the Leray spectral sequence for the fibration $\pi_2\co  {\I_2}\rightarrow \F2{\Pp2}$. It is simpler to consider the $\C^*$-quotient and study the fibration $\pi_2'\co  \mathcal P_2\rightarrow \F2{\Pp2}$, which is a locally trivial fibration with fibre $\Pp{}(V_p\cap V_q\setminus\Sigma)$. The $E_2$ terms of the Leray spectral sequence are written in Table~\ref{lera}. 
Note that the space $\F2{\Pp2}$ is simply connected and that its cohomology has $\s_2$-HG polynomial $(\inv+\Ll t\ant)(1+\Ll t^2+\Ll^2t^4)$ with respect to the natural action of $\s_2$ generated by the involution $(\alpha,\beta)\leftrightarrow(\beta,\alpha)$. 

\begin{table}\caption{\label{lera} First block of the Leray spectral sequence in cohomology associated to $\pi_2'$}
$$\begin{array}{r|@{}ccccccc@{}}
7 &\begin{smallmatrix}\Ant\\\otimes\\\Q(-5)\end{smallmatrix}
&0&\begin{smallmatrix}(\Inv+\Ant)\\\otimes\\\Q(-6)\end{smallmatrix}
&0&\begin{smallmatrix}(\Inv+\Ant)\\\otimes\\\Q(-7)\end{smallmatrix}
&0&\begin{smallmatrix}\Inv\\\otimes\\\Q(-8)\end{smallmatrix}
\\[12pt]
6 &\begin{smallmatrix}(\Inv+\Ant)\\\otimes\\\Q(-4)\end{smallmatrix}
&0&\begin{smallmatrix}(\bigoplus^2\Inv+\bigoplus^2\Ant)\\\otimes\\\Q(-5)\end{smallmatrix}
&0&\begin{smallmatrix}(\bigoplus^2\Inv+\bigoplus^2\Ant)\\\otimes\\\Q(-6)\end{smallmatrix}
&0&\begin{smallmatrix}(\Inv+\Ant)\\\otimes\\\Q(-7)\end{smallmatrix}
\\[17pt]
5 &0&0&0&0&0&0&0\\[12pt]
4 &\begin{smallmatrix}\Ant\\\otimes\\\Q(-3)\end{smallmatrix}
&0&\begin{smallmatrix}(\Inv+\Ant)\\\otimes\\\Q(-4)\end{smallmatrix}
&0&\begin{smallmatrix}(\Inv+\Ant)\\\otimes\\\Q(-5)\end{smallmatrix}
&0&\begin{smallmatrix}\Inv\\\otimes\\\Q(-6)\end{smallmatrix}
\\[12pt]
3 &\begin{smallmatrix}(\bigoplus^2\Inv+\Ant)\\\otimes\\\Q(-2)\end{smallmatrix}
&0&\begin{smallmatrix}(\bigoplus^3\Inv+\bigoplus^3\Ant)\\\otimes\\\Q(-3)\end{smallmatrix}
&0&\begin{smallmatrix}(\bigoplus^3\Inv+\bigoplus^3\Ant)\\\otimes\\\Q(-4)\end{smallmatrix}
&0&\begin{smallmatrix}(\Inv+\bigoplus^2\Ant)\\\otimes\\\Q(-5)\end{smallmatrix}
\\[17pt]
2 &0&0&0&0&0&0&0
\\[17pt]
1 &0&0&0&0&0&0&0\\[12pt]
0 &\Inv&0&\begin{smallmatrix}(\Inv+\Ant)\\\otimes\\\Q(-1)\end{smallmatrix}&0&\begin{smallmatrix}(\Inv+\Ant)\\\otimes\\\Q(-2)\end{smallmatrix}&0&\begin{smallmatrix}\Ant\\\otimes\\\Q(-3)\end{smallmatrix}\\\hline
&0&1&2&3&4&5&6
\end{array}$$
\end{table}

The proof of the following lemma is based on a suggestion by Alexei Gorinov.

\begin{lem}\label{lempi} In the spectral sequence associated to $\pi_2'$, the differential $d_4\co  E_4^{0,3}\rightarrow E_4^{4,0}$ has rank two.
\end{lem}

\begin{proof}
Denote by $\mathcal P_1\subset \Pp2\times\Pp{}(S^2_4\setminus\Sigma)$ the variety of pairs $(\xi,[f])$ such that $f(\xi)=0$. Consider the inclusion $i\co  \mathcal P_2\rightarrow\mathcal P_1\times\mathcal P_1$ defined by $i(\alpha,\beta,[f])=((\alpha,[f]),(\beta,[f]))$. There is a commutative diagram
\begin{equation}\label{pi}
\begin{CD} \mathcal P_2 @>i>> \mathcal P_1\times\mathcal P_1\\
@V\pi_2'VV @VV\pi_1\times\pi_1V\\
\F2{\Pp2} @>\text{inclusion}>> \Pp2\times\Pp2,
\end{CD}\end{equation}
where $\pi_1\co  \mathcal P_1\rightarrow\Pp2$ denotes the natural projection.

In particular, the differentials of the spectral sequences associated to $\pi_2'$ and $\pi_1\times\pi_1$ commute with the maps induced by $i$ on the $E_4$ terms of the spectral sequences. 

Recall from \cite[Proposition~1]{BT} that the differential 
$$E^{0,3}_4(\pi_1)\xrightarrow{d_4} E^{4,0}_4(\pi_1)$$
has rank $1$. This implies that the differential $E^{0,3}_4(\pi_1\times\pi_1)\xrightarrow{d_4} E^{4,0}_4(\pi_1\times\pi_1)$ has rank~2. Since 
$$E_4^{4,0}(\pi_1\times\pi_1)\cong E_2^{4,0}(\pi_1\times\pi_1)\cong \coh4{\Pp2\times\Pp2}\otimes\coh0{V_p\setminus\Sigma}$$
and 
$$E_4^{4,0}(\pi_2')\cong E_2^{4,0}(\pi_2')\cong \coh4{\F2{\Pp2}}\otimes\coh0{V_p\cap V_q\setminus\Sigma},$$
one can verify directly that the composition of $d_4\co  E^{0,3}_4(\pi_1\times\pi_1)\rightarrow E^{4,0}_4(\pi_1\times\pi_1)$ and the map $E^{4,0}_4(\pi_1\times\pi_1)\rightarrow E^{4,0}_4(\pi_2')$ is surjective. Then the claim follows from the commutativity of the diagram~\ref{pi}.
\end{proof}

\begin{proof}[Proof of Theorem~\ref{result}]
Comparing the Leray--Hirsch isomorphism~\eqref{icLH} with Table~\ref{lera} implies that the entire contribution of the first block to the cohomology of $\I_2$ is determined by the cohomology of $\I_2$ in degree $\leq5$. This follows from the fact that the cohomology of $\GL(3)$ is trivial in degree $k\geq10$. 
Then Lemma~\ref{lempi}, together with the structure of $\coh\pu{\I_2}$ as tensor product of $\GL(3)$, yields that the first block contributes 
$$(1+\Ll t)(1+\Ll^2 t^3)(1+\Ll^3 t^5)(\inv+\Ll t^2(\inv+\ant)+\Ll^3 t^5 \inv).$$to the $\s_2$-HG polynomial of $\coh\pu{\I_2}$.
This implies that the $\s_2$-HG polynomial of the cohomology of the moduli space $\mathcal Q_2$ of smooth quartic curves with two marked points is 
\begin{equation}\label{firstblock}
\inv+\Ll t^2(\inv+\ant)+\Ll^3 t^5 \inv,
\end{equation}
plus the term coming from singular configurations of type $8$--$10$. In the next section (see page \pageref{secondblock}), we will prove that this term equals
\begin{equation*}\tag{\ref{secondblock}}
\Ll^6t^6\inv+\Ll^7t^8(\inv+\ant)+\Ll^8t^8\ant.\end{equation*}

Summing the contributions of the two block of columns, we get that the $\s_2$-HG polynomial of the cohomology of $\mathcal Q_2$ is
$$\inv+\Ll t^2(\inv+\ant)+\Ll^3 t^5 \inv+\Ll^6t^6\inv+\Ll^7t^8(\inv+\ant)+\Ll^8t^8\ant.$$

This establishes the first part of Theorem~\ref{result}. 
To prove the second part of the theorem, recall from~\cite[Corollary III.2.2]{OTtesi} that the $\s_2$-HG polynomial of the cohomology of the hyperelliptic locus $\Hhm32\subset\Mm32$ is $\inv+\Ll t^2(\inv+\ant)+ \Ll^7t^7\ant$, and consider the long exact sequence associated to the inclusion $\Hhm32\hookrightarrow\Mm32$:
$$\cdots\rightarrow\coh k{\Mm32}\rightarrow\coh k{\mathcal Q_2}\rightarrow\coh{k-1}{\Hhm32}\otimes\Q(-1)\rightarrow\coh{k+1}{\Mm32}\rightarrow\cdots,$$
which can be rephrased in Borel--Moore homology as 
$$\cdots\rightarrow\bar H_k(\Mm32;\Q)\rightarrow\bar H_k(\mathcal Q_2;\Q)\xrightarrow{d_k}\bar H_{k-1}(\Hhm32;\Q)\rightarrow\bar H_{k-1}(\Mm32;\Q)\rightarrow\cdots.$$

If $k\neq8$, the differentials $d_k$ are always zero for Hodge-theoretic reasons.  If $k=8$, both $\bar H_8(\mathcal Q_2;\Q)$ and $\bar H_7(\Hhm32;\Q)$ have a one-dimensional summand of Hodge weight $0$, on which $\s_2$ acts as the sign representation. Hence, a priori $d_8$ can have either rank $0$ or $1$. 
To determine the rank of $d_8$, we observe that both the cohomology of $\Hhm32$ and $\mathcal Q_2$ were computed using Vassiliev--Gorinov's method. In particular, both $\bar H_8(\mathcal Q_2;\Q)$ and $\bar H_7(\Hhm32;\Q)$ are related to configurations of at least $4$ singular points.
Moreover, there configurations correspond to strata of the geometric realizations that have Borel--Moore homology which is a tensor product of that of the group acting. This means that both Borel--Moore homology groups can be interpreted as Borel--Moore homology groups of certain moduli spaces.

Specifically, consider the moduli space $\mathcal N$ whose elements are isomorphism classes of triples $(C,p,q)$, where $C$ is the union of two smooth rational curves intersecting transversally at $4$ distinct points and $p,q$ are any distinct (but possibly singular) points on $C$. Note that the arithmetic genus of such a curve $C$ is $3$. Denote by $S$ the rank $1$ local system on $\mathcal N$ changing its orientation every time a pair of nodes on $C$ is interchanged, and denote by $\mathcal N_h$ the closed subset of $\mathcal N$ such that the four nodes have the same moduli on both rational components. 

Observe that the problem with the determination of $d_8$ only concerns the Hodge weight 0 summands of the Borel-Moore homology groups. For this reason, in the rest of the proof we will restrict to the Hodge weight 0 summands of each homology group we consider. 

The space $\mathcal N$ can be written as the disjoint union of locally closed strata, each of them isomorphic to the quotient by the action of a finite group, of a product of moduli spaces $\Mm0n$ with $4\leq n\leq 6$, whose cohomology groups are completely known (see e.g. \cite{G-0}). Investigating this stratification, one gets that the only Borel--Moore homology group of Hodge weight 0 is $\bar H_4(\mathcal N;S)=\Ant$. Analogous considerations also apply to $\mathcal N_h$. In that case, one has $\bar H_3(\mathcal N_h;S)\cong\Ant$ as only Borel--Moore homology group with Hodge weight $0$. By the constructions in~\cite[III.2]{OTtesi}, there is a natural isomorphism $\bar H_3(\mathcal N_h;S)\cong\bar H_7(\Hhm 32;\Q)$.

The weight $0$ part of the Borel--Moore homology of $\mathcal N\setminus\mathcal N_h$ can also be computed directly with Vassiliev--Gorinov's method. This yields again that the only non-trivial Borel--More homology group with $S$-coefficient of $\mathcal N\setminus\mathcal N_h$ is $\bigoplus_2\Ant$ in degree $4$. Moreover, the direct computation shows that $\bar H_4(\mathcal N\setminus\mathcal N_h;S)$ is generated by two classes, both related to configurations of type (9) in Table~\ref{talist}. This allows to define a surjective map $\bar H_4(\mathcal N\setminus\mathcal N_h;S)\rightarrow\bar H_8(\mathcal Q_2;\Q)$ making the following diagram commute:
$$
\begin{CD}
0@>>>{\stackrel{\begin{array}{c}\Ant\\\|\wr\end{array}}{\bar H_4(\mathcal N;S)}}
@>>>{\stackrel{\begin{array}{c}\bigoplus_2\Ant\\\|\wr\end{array}}{\bar H_4(\mathcal N\setminus\mathcal N_h;S)}}
@>{\gamma}>>{\stackrel{\begin{array}{c}\Ant\\\|\wr\end{array}}{\bar H_3(\mathcal N_h;S)}}@>>>0\\
@.@.@VVV@VVV\\
@.@.\bar H_8(\mathcal Q_2;\Q)@>{d_8}>>\bar H_7(\Hhm 32;\Q).
\end{CD}
$$
The commutativity of the diagram immediately yields $\rank d_8=\rank \gamma=1$. 
\end{proof}

\section{Configurations of five and six points}\label{QKQK}

The aim of this section is to compute the contribution of singular configurations of type $8$, $9$ and $10$ (see Table~\ref{talist}) to the rational cohomology of $\I_2$ and $\mathcal Q_2$. 
For these configuration, it seems more natural to work directly with the cohomology of $\I_2$, without having to pass through the study of the fibre of $\pi_2$. This is indeed possible.
Namely, consider the space
$$\mathcal D:=\{(\alpha,\beta,f)\in\F2{\Pp2}\times\Sigma: f(\alpha)=f(\beta)=0 \}.$$

Note that $\mathcal D$ is a closed subset of $\mathcal V:=\{(\alpha,\beta,f)\in\F2{\Pp2}:f(\alpha)=f(\beta)=0\}$. The space $\mathcal V$ is the total space of a vector bundle over $\Pp2$, and $\I_2=\mathcal V\setminus\mathcal D$. 
Vassiliev--Gorinov's method can be exploited to compute the Borel--Moore homology of $\mathcal D$. This is done by defining the singular locus of an element $(\alpha,\beta,f)$ in $\mathcal D$ as the subset $\{(\alpha,\beta)\}\times K_f$ of $\F2{\Pp2}\times \Pp2$, where $K_f$ denotes the singular locus of the polynomial $f$. In particular, the classification of singular sets of elements of $\mathcal D$ is obtained from the classification of singular sets of elements of $V_p\cap V_q\setminus\Sigma$ by allowing the pair $(p,q)$ to move in $\F2{\Pp2}$.

Even though this is no longer the original setting of Vassiliev--Gorinov's method, one can mimic the construction of the cubical spaces $\Lambda$ and $\mathcal X$ (see \S~\ref{VGmethod}), and obtain cubical spaces $\Lambda'$ and $\mathcal X'$ that play an analogous role. In particular, the map $\ba\mathcal X'\ba\rightarrow \mathcal D$ induces an isomorphism on the Borel--Moore homology of these spaces, because it is a proper map with contractible fibres.
Moreover, for the stratifications $\Phi'$ and $F'$ obtained from the construction of $\Lambda'$ and $\mathcal X'$, we have natural maps $\Phi'_k\rightarrow\F2{\Pp2}$ and $F'_k\rightarrow\F2{\Pp2}$ which are locally trivial fibrations with fibre isomorphic to $\Phi_k$, respectively, $F_k$. 

In view of the considerations above, computing the rational Borel--Moore homology of the spaces $\Phi'_j$ and $F'_j$ for $j\in\{8,9,10\}$ is enough to get the contribution of configurations of type $8$--$10$ to the cohomology of $\I_2$.
We start by determining the twisted Borel--Moore homology of the underlying families of configurations $X'_8$ and $X'_9$.

Define $Y_8\subset\F2{\Pp2}\times\F5{\Pp2}$ 
to be the space of configurations $(p,q,e_1,e_2,e_3,$ $e_4,c)$ such that \begin{itemize}
\item $\{c\}=e_1e_2\cap e_3e_4\not\subset\{p,q\}$;
\item $p,q,e_1,e_2,e_3,e_4$ lie on a conic different from the reducible conic $e_1e_2\cup e_3e_4$;
\item $\{p,q\}\not\subset\{e_1,e_2,e_3,e_4\}$.
\end{itemize}
Then $X'_8\subset\F2{\Pp2}\times\B5{\Pp2}$ is isomorphic to the quotient of $Y_8$ by the action of the subgroup $G$ of $\s_4$ generated by the permutations $(1,2)$, $(3,4)$ and $(1,3)(2,4)$. The action of $G$ on $Y_8$ is given by permuting the four points $(e_1,e_2,e_3,e_4)$ in the configurations.
Since $G$ is a subgroup of $\s_4$, it makes sense to restrict the sign representation to it.

Furthermore, note that the conic passing through the points $p,q,e_1,e_2,e_3,e_4$ is uniquely determined for every configuration in $Y_8$. Therefore, $Y_8$ can be embedded in the space $W\subset\F2{\Pp2}\times\F5{\Pp2}\times\Pp{}(S^2_2)$ of configurations  $(p,q,e_1,e_2,e_3,e_4,c,$ $C)$, such that 
\begin{itemize}
\item $\{c\}=e_1e_2\cap e_3e_4\not\subset\{p,q\}$;
\item the points $p,q,e_1,e_2,e_3,e_4$ lie on the conic $C$;
\item the conic $C$ is distinct from the reducible conic $e_1e_2\cup e_3e_4$.
\end{itemize}
Hence, we have the chain of inclusions
$$
Y_8\hookrightarrow W\hookrightarrow \F2{\Pp2}\times\F5{\Pp2}\times\Pp{}(S_2^2).$$

\begin{lem}\label{X8}
Denote by $S$ the local system of coefficients induced on $W/G$ and $X'_8$ by the sign representation on $G$.
Consider the $\s_2$-action generated by the involution interchanging the points $(p,q)\in\F2{\Pp2}$.
Then one has
$$\Phgs(\bar H_\pu(X'_8;S))=(\Ll^{-1} t^3+\Ll^{-2} t^4)\inv\cdot\Phg(\bar H_\pu(\PGL(3);\Q)).$$
\end{lem}

To prove Lemma~\ref{X8}, we will consider the quotient of $Y_8$ and $W$ by the action of $\PGL(3)$. Since every configuration in $W$ contains points $e_1,e_2,e_3,e_4$ which are in general position, the group $\PGL(3)$ acts freely and transitively on $W$, hence $W$ is isomorphic to the product of $\PGL(3)$ and the quotient $W/\PGL(3)$. Recall that $\PGL(3)$ is isomorphic to the configuration space of four ordered points in general position in $\Pp2$. This yields a natural identification between the quotient $W/\PGL(3)$ and the space 
$$W_E:=W\cap(\F2{\Pp2}\times\{(E_1,E_2,E_3,E_4,E_5)\}),$$
where 
$$E_1=[1,0,0],\ \ E_2=[0,1,0],\ \ E_3=[0,0,1],\ \ E_4=[1,1,1],\ \ E_5=[1,1,0].$$
In the following, we identify each element $\sigma$ of $G$ with
the automorphism of $\Pp2$ mapping $E_i$ to $E_{\sigma(i)}$. 
This allows to consider $G$ as a subgroup of $\mathrm{Aut}(\Pp2)$, and induces an action of $G$ on $W_E$ that makes the isomorphism $W\cong W_E\times\PGL(3)$ is $G$-equivariant. 
Note that the action of $G$ on $\PGL(3)$ is defined by restricting to $G$ the natural action of the symmetric group $\s_4$ on $\PGL(3)\hookrightarrow\F4{\Pp2}$ permuting the four points in the configuration. It is not difficult to prove that the rational Borel--Moore homology of $\PGL(3)$ is $\s_4$-invariant and hence also $G$-invariant. 

By applying the K\"unneth formula to the Borel--Moore homology of $W\cong W_E\times\PGL(3),$ and considering the part of the Borel--Moore homology which has the wished behaviour for the action of $G$, one gets
\begin{equation}\label{equaW}
\bar H_\pu(W/G;S)\cong \bar H_\pu(W_E/G;S)\otimes\bar H_\pu(\PGL(3);\Q),
\end{equation}
where $S$ denotes the local system of rank $1$ induced by the restriction of the sign representation to $G\subset\s_4$. 

The reasoning above applies to $Y_8$ as well as $W$. This yields the isomorphism
\begin{equation}\label{equaY8}
\bar H_\pu(X'_8;S)\cong \bar H_\pu((Y_8\cap W_E)/G;S)\times\bar H_\pu(\PGL(3);\Q).
\end{equation}

The space $W_E$ can be described in the following way. 
Denote by $\mathcal L$ the space of conics passing through the $E_i$'s and distinct from the reducible conic $(x_0-x_1)x_2=0$. Note that $\mathcal L$ is isomorphic to an affine line. Then we have
$$W_E=\{(p,q,C)\in\F2{\Pp2}\times \mathcal L: p,q\in C\}.$$

\begin{lem}\label{utile}
In the notation of Lemma~\ref{X8}, we have $\Phgs(\bar H_\pu(W_E/G;S))=\Ll^{-2}t^4\inv.$
\end{lem}

\begin{proof}
The space $Q=W_E/G$ can be decomposed as the union of a closed locus $K$ containing all equivalence classes of triples $(p,q,C)$ such that $C$ is a singular conic, and an open part $U$ where the conic $C$ is always non-singular. 

We compute the Borel--Moore homology of $K$ first. The locus $K$ has two components, according to the position of the two points $p,q$. We denote by $M$ the component of $K$ such that $p,q$ lie on the same irreducible component of $C$, and $N$ the component in which $p,q$ lie on two different components of $C$. The elements of the intersection $M\cap N$ are the configurations in which the singular point of $C$ is either $p$ or $q$.

Up to the $G$-action, the space $M$ can be identified with the space of ordered configurations of two points on the projective line $x_1=0$, hence the $\s_2$-HG polynomial of $\bar H_\pu(M;S)$ is $(\Ll^{-1}t^2+\Ll^{-2}t^4)\inv$.

Next, we identify the space $N\setminus M$ with a $\s_2$-quotient of the space of pairs $(p,q)$ where $p$ lies on $x_1=0$, the point $q$ lie on $x_0-x_2=0$ and both points are distinct from the intersection point of these lines. The $\s_2$-action interchanges $E_1$ and $E_3$, and $E_2$ and $E_4$, and we have to take invariant classes with respect to it. This implies that the $\s_2$-HG polynomial of $\bar H_\pu(N\setminus M;S)$ is $(\Ll^{-2}t^4)\inv$. 
Then, from the long exact sequence in Borel--Moore homology associated with the closed inclusion $M\hookrightarrow K$ we can conclude that the $\s_2$-HG polynomial of $\bar H_\pu(K;S)$ is $(\Ll^{-1}t^2+2\Ll^{-2}t^4)\inv$.

Subsequently, we compute the Borel--Moore homology of $U$ by lifting $U$ to a $G$-invariant subset $U'\subset W_E$, and looking for the part of the Borel--Moore homology of $U'$ that has the wished behaviour with respect to the $G$-action. We have that $U'$ projects to the locus of non-singular conics in $\mathcal L$, which is isomorphic to $\C\setminus\{\pm1\}$. Note that the action of $(1,2)\in G$ on $\mathcal L\cong \C$ interchanges the two singular conics. The projection $U'\rightarrow\C\setminus\{\pm1\}$ is a locally trivial fibration with fibre isomorphic to the space $\F2R$, where $R$ is a chosen non-singular conic through the $E_i$'s. In order to study the action of $G$ on the Borel--Moore homology of $\F2R$, we assume that the conic $R$ is fixed by all automorphisms in $G\subset\Aut(\Pp2)$. If we fix an isomorphism $R\cong\Pp1$, we have that taking the quotient by $G$ gives finite maps $R\cong\Pp1\rightarrow\Pp1$ and $\F2R\cong\F2{\Pp1}\rightarrow\F2{\Pp1}$. In particular, the Borel--Moore homology with standard coefficients of $\F2R$ is isomorphic to that of its quotient by $G$, hence all Borel--Moore homology classes of $\F2R$ are $G$-invariant.
Hence, the $\s_2$-HG polynomial of $\bar H_\pu(U;S)$ is the product of $\Phg(\F2{\Pp1})\inv$ and the HG polynomial of the part of the Borel--Moore homology of $\C\setminus\{\pm1\}$ which is anti-invariant for the involution $\xi\leftrightarrow -\xi$, which equals $t$. 

To compute $\bar H_\pu(Q;S)$, we can now use the long exact sequence in Borel--Moore homology associated to the closed inclusion $K\hookrightarrow Q$:
$$\cdots\rightarrow\bar H_k(K;S)\rightarrow\bar H_k(Q;S)\rightarrow\bar H_k(U;S)\rightarrow\bar H_{k-1}(K;S)\rightarrow\cdots$$

This yields immediately $\bar H_k(Q;S)=0$ if $k>5$ or $k<2$. Moreover, we have
$$
0\rightarrow\bar H_5(Q;S)\rightarrow\Q(2)\xrightarrow{\delta_5}\Q(2)^2\xrightarrow{\delta_5'}\bar H_4(Q;S)\rightarrow 0,$$
$$0\rightarrow\bar H_3(Q;S)\rightarrow\Q(1)\xrightarrow{\delta_3}\Q(1)\xrightarrow{\delta_3'}\bar H_2(Q;S)\rightarrow 0.$$

Then the claim follows from the fact that both $\delta_5$ and $\delta_3$ are injections. As we will see, the subset $M\cup U$ has trivial Borel--Moore homology with $S$-coefficients, hence $\bar H_k(M;S)$ is contained in the kernel of $\delta_k'$ for every $k$.

To compute the Borel--Moore homology of $M\cup U$, consider the surjective map $\pi\co M\cup U\rightarrow\mathcal L/G$ obtained by restricting the natural projection $W_E\rightarrow \mathcal L$. 
The map $\pi$ is clearly locally trivial on $U$. We claim that $\pi$ is also locally trivial in a neighborhood of the point $w_0$ in $\mathcal L/G$ parametrizing singular conics. Up to the $G$-action, and possibly the choice of a sufficiently small neighbourhood $U_0$, we can assume that this singular conic is $Y\co  x_1(x_0-x_2)=0$, and identify $\pi^{-1}(w)$ ($w\in U_0$) with the locus of triples $(Y,\alpha,\beta)$ such that $\alpha$ and $\beta$ lie on the line $x_1=0$. Then the fibre of $\pi$ near $w_0$ can be identified with $\F2{\{x_1=0\}}$ by considering the projection from the point $E_4$, which maps every non-singular conic in $\mathcal L$ onto the line ${x_1=0}$. This construction yields a map from the preimage in $\pi$ of a neighbourhood of $w_0$ to $\F2{\Pp1}\cong\F2{\{x_1=0\}}$, which admits a section. Hence, the map $\pi$ is a locally trivial fibration over $\mathcal L/G$.

Note that this implies that the Borel--Moore homology of $M\cup U$ in the local system $S$ is trivial.
The Borel--Moore homology of $\mathcal L$ is clearly $G$-invariant, hence the elements of $H_\pu(M\cup U;S)$ have to come from Borel--Moore homology classes of the fibre of $\pi$, in a local system different from the standard one.
The fibre of $\pi$ is isomorphic to a $G$-quotient of $\F 2{\Pp1}$, and the whole Borel--Moore homology of $\F2{\Pp1}$ is $G$-invariant. For this reason, the fibre of $\pi$ has trivial Borel--Moore homology in all local systems different from the standard one.
\end{proof}

\begin{proof}[Proof of Lemma \ref{X8}]

We start by investigating the space $W_E\setminus Y_8$ and its quotient $Q'$ by the action of $G$. Recall that a configuration $(p,q,C)\in Q$ lies in $Q'$ if and only if $\{p,q\}\subset\{E_1,E_2,E_3,E_4\}$. It is easy to see that $Q'$ has two components, according to whether $p$ and $q$ lie both on the same component of the reducible conic $x_2(x_0-x_1)=0$, or not. Denote by $Q_a$ the component corresponding to the first case and by $Q_b$ the component corresponding to the second case. Up to the action of $G$, we may assume that for every configuration in $Q_a$ we have $p=E_1$, $q=E_2$. Hence, the space $Q_a$ is isomorphic to the quotient $\mathcal L/\iota$, where the involution $\iota$ is $(3,4)\in G\subset \Aut(\Pp2)$. Since $\mathcal L/\iota$ is isomorphic to $\C$, the Borel--Moore homology of $Q_a$ with $S$-coefficients is isomorphic to the Borel--Moore homology of $\C$ induced by the sign representation on $\langle\iota\rangle = \s_2$, which is trivial. 

Analogously, up to the $G$-action one can assume that $p=E_1$, $q=E_3$ hold for every configuration in $Q_b$. In particular, $Q_b$ is isomorphic to $\mathcal L\cong \C$ and is invariant for the involution interchanging $p$ and $q$. This proves that $\bar H_\pu(Q';S)$ is isomorphic to $\bar H_\pu(\C)$ and is invariant for the involution $p\leftrightarrow q$. Then the claim follows from the long exact sequence in Borel--Moore homology associated to the closed inclusion $Q'\hookrightarrow Q$ and isomorphism~\eqref{equaY8}. 
\end{proof}

Recall from Table~\ref{talist} that $X'_9$ is the locus 
$$X'_9:=\big\{(p,q,S)\in \F2{\Pp2}\times\B6{\Pp2}:\exists \{r_i\}_{1\leq i\leq4}\in\B4{\Pp2\duale}\big(S=Sing(\bigcup_ir_i), p,q\in\bigcup_ir_i\big)\big\}.
$$

Observe that giving six points that are the pairwise intersection of four lines in general position is equivalent to giving the configuration of four lines.
Denote by $\tF4{\Pp2\duale}$ the space of ordered configurations of lines in general position (i.e., such that no three of them pass through the same point), and by $\tB4{\Pp2\duale}$ the analogous space of unordered configurations. Then we have
$$
X'_9\cong\left\{(p,q,\{r_1,r_2,r_3,r_4\})\in \F2{\Pp2}\times\tB4{\Pp2\duale}:  p,q\in\bigcup_ir_i\right\}.
$$

We start by investigating the closed subset $X'_{9a}$ of configurations $(p,q,\{r_i\}_i)\in X'_9$ such that $p$ and $q$ lie on the same line $r_j$ for some index $j$. 

\begin{lem}\label{X9}
$$\Phgs(\bar H_\pu(X'_{9a};\pm\Q))=\Phgs(\bar H_\pu(\F2{\Pp1};\Q))\cdot\Phg(\bar H_\pu(\PGL(3);\Q)).$$
\end{lem}

\begin{proof}
Consider the variety
$$A:=\{(p,q,r_1,r_2,r_3,r_4)\in \F2{\Pp2}\times\tF4{\Pp2\duale}: p,q\in r_4\}.$$
Note that $X'_{9a}$ is the quotient of $A$ by the action of $\s_3$ interchanging $r_1, r_2$ and $r_3$. On the other hand, we have $A\cong\F2{\Pp1}\times\PGL(3)$, where we used the fact that $\PGL(3)$ is isomorphic to the space of four lines in general position, and chosen an isomorphism $\Pp1\cong r_4$ (for instance, the one mapping $0$ to $r_1\cap r_4$, $1$ to $r_2\cap r_4$ and $\infty$ to $r_3\cap r_4$). Hence, we can obtain the Borel--Moore homology of $X'_{9a}$ by taking the $\s_3$-invariant part of the Borel--Moore homology of $A$. This establishes the claim.
\end{proof}

Next, we consider $X'_{9b}:=X'_9\setminus X'_{9a}$.

\begin{lem}\label{X10}
$$\Phgs(\bar H_\pu(X'_{9b};\pm\Q))=(t^2\ant+\Ll^{-2}t^4\inv)\cdot\Phg(\bar H_\pu(\PGL(3);\Q)).$$
\end{lem}

\begin{proof}
Consider the space 
$$Y_9:=\{(p,q,r_1,r_2,r_3,r_4)\in\F2{\Pp2}\times\tF4{\Pp2\duale}:p,q\in\bigcup_{i}r_i\}.$$

Observe that $\tF4{\Pp2\duale}$ is isomorphic to $\PGL(3)$, and that the group $\PGL(3)$ acts freely and transitively on $Y_9$. 
The quotient of this action is isomorphic to the fibre of the projection $Y_9\rightarrow\tF4{\Pp2\duale}$ at the configuration $(l_1,l_2,l_3,l_4)$, where
$$l_1\co  x_0=0,\ \ l_2\co  x_1=0,\ \ l_3\co  x_2=0,\ \ l_4\co  x_0+x_1+x_2=0.$$If we pose $L\co  x_0x_1x_2(x_0+x_1+x_2)=0$, this implies that $Y_9/\PGL(3)$ is isomorphic to $\F2L$, and we have an isomorphism 
\begin{equation}\label{isoy}
Y_9\cong\F2L\times\PGL(3).
\end{equation}

Consider the action of $\s_4$ on $L$ and $\F2L$ defined by identifying every permutation $\sigma\in\s_4$ with the automorphism of $\Pp2$ sending the line $l_i$ to $l_{\sigma_i}$ for all $i$, $1\leq i\leq4$. The natural action of $\s_4$ on $\tF4{\Pp2}$ defines an action on $\PGL(3)$ via the isomorphism $\PGL(3)\cong\tF4{\Pp2\duale}$, making isomorphism~\eqref{isoy} $\s_4$-equivariant. 
Applying K\"unneth formula and taking the $\s_4$-invariant part of the Borel--Moore homology of $Y_9$ yields
$$\bar H_\pu(X'_{9};\Q)\cong \bar H_\pu(\F2L;\Q)^{\s_4}\otimes\bar H_\pu(\PGL(3);\Q),$$
where we used the fact that $X'_9$ is the quotient of $Y_9$ by the action of $\s_4$, and that the whole Borel--Moore homology of $\PGL(3)$ is $\s_4$-invariant.

Let us see what these considerations tell us about the Borel--Moore homology of $X'_{9b}$. The space $X'_{9b}$ is isomorphic to the $\s_4$-quotient of the product of $\PGL(3)$ and the locus of configurations of two points $(a,b)\in\F2L$ not lying on the same component of $L$. This locus can be decomposed according to whether $a$ and $b$ are or are not singular points of $L$ into the loci
$$S_1:=\{(a,b)\in\F2L:\text{$a$ and $b$ are both singular points}\},$$
$$S_2:=\{(a,b)\in\F2L:\text{only one of the points $a$ and $b$ is a singular point of }L\},$$
$$S_3:=\{(a,b)\in\F2L:\text{$a$ and $b$ are non-singular points of }L\}.$$

The quotient $S_1/\s_4$ consists of only one point, the class of the pair $([1,0,0],$ $[0,1,-1])$. The quotient $S_2/\s_4$ has two isomorphic components, according to which point ($a$ or $b$) is a singular point of $L$. Consider the case in which $a$ is singular. Up to the action of $\s_4$, we can assume that $a$ is the point $[1,0,0]$ and $b$ lies on $x_0=0$. By the definition of $S_2$ we know that $b$ is different from the points $[0,1,-1]$, $[0,0,1]$ and $[0,1,0]$. Note that, since we are working modulo $\s_4$, the coordinates of $b$ are defined up to the involution interchanging $x_1$ and $x_2$. This proves that both components of $S_2/\s_4$ are isomorphic to $\C^*$. 

Finally, we determine the Borel--Moore homology of the quotient of $S_3$ by the action of $\s_4$. Up to the action of the group, we can assume that $a$ lies on the line $l_3$ and $b$ on $l_4$. The position of both points is determined up to the involution interchanging the lines $l_1$ and $l_2$. If we identify $l_3$ and $l_4$ with $\Pp1$, and $l_3\cap l_4$ with the point at infinity of the projective line, we have that $S_3/\s_4$ can be embedded into the quotient of $(\C\setminus\{\pm1\})^2$ by the relation $(t,s)\sim(-t,-s)$. The complement of $S_3/\s_4$ in this quotient is the locus such that either $t$ or $s$ are equal to $\pm1$. We can study $(\C\setminus\{\pm1\})^2/\sim$  as follows:
$$\begin{array}{ccccc}
\C^2 & \xrightarrow{\mod \sim} & \{(x,y,z)\in\C^3: y^2=xz\} & \xrightarrow{\mod \s_2}&\C^2\\
(t,s)&\longmapsto & (t^2,ts,s^2) & \longmapsto &(t^2+s^2,ts)\\
(1,s)&\longmapsto & (1,s,s^2) & \longmapsto &(s^2+1,s)\\
(t,1)&\longmapsto & (t^2,t,1) & \longmapsto &(t^2+1,t),
\end{array}$$
where the second map denotes the quotient by the action of $\s_2$ interchanging $t$ and $s$.
Concluding, the spectral sequence associated to this stratification has $E^1$ term  as in Table~\ref{littlespec} (where we have taken into account the $\s_2$-action interchanging $a$ and $b$).

\begin{table}\caption{\label{littlespec}}
$$\begin{array}{r|ccc}
 1&0&0&\Inv\otimes\Q(2)\\
 0&0&(\Inv+\Ant)\otimes\Q(1)&(\Inv+\Ant)\otimes\Q(1)\\
-1&\Inv&\Inv+\Ant&\bigoplus^2\Ant\\\hline
&1&2&3
\end{array}
$$
\end{table}

We can use the geometric description of $S_1$, $S_2$ and $S_3$ to determine all differentials of the spectral sequence above.
In particular, the the $0$-th row is exact, and both differentials in the row of index $-1$ have rank 1. Then the claim follows from the fact that $\bar H_\pu(X'_{9b};\Q)$ is isomorphic to  $\bar H_\pu(S_1\cup S_2\cup S_3;\Q)^{\s_4}\otimes\bar H_\pu(\PGL(3))$.
\end{proof}

\begin{prop}
The Borel--Moore homology groups of the unions of strata $\Phi'_8\cup \Phi'_9\subset\ba\Lambda'\ba$ and $F'_8\cup F'_9\cup F'_{10}\subset\ba\X'\ba$ are as follows:
\begin{equation}\label{een}
\Phgs(\bar H_\pu(\Phi'_8\cup \Phi'_9;\Q))=(\Ll^{-2}t^9\inv+(\Ll^{-1}\inv+\Ll^{-1}\ant+\ant)t^7)\cdot\Phg(\bar H_\pu(\PGL(3);\Q)).\end{equation}
\begin{equation}
\label{twee}
\Phgs(\bar H_\pu(F'_8\cup F'_9\cup F'_{10};\Q))=(\Ll^{-2}t^9\inv +(\Ll^{-1}\inv+\Ll^{-1}\ant+\ant)t^7)\cdot\Phg(\bar H_\pu(\GL(3);\Q)).\end{equation}
\end{prop}

\begin{proof}
Lemmas~\ref{X9} and~\ref{X10} imply that $\Phgs(\bar H_\pu(X'_9;\pm\Q))$ equals $(1+\Ll^2t^{-3})(1+\Ll^3t^{-5})\Ll^{-8}t^{16}\cdot(2\Ll^{-2}t^4\inv +(\Ll^{-1}+1)t^2\ant)$. Recall from Table~\ref{talist} (page~\pageref{talist}) that $\Phi'_9$ is a simplices bundle with $5$-dimensional fibre. Hence, the $\s_2$-HG polynomial of the Borel--Moore homology of $\Phi'_9$ equals that of $X'_9$ multiplied by $t^5$. Analogously, Lemma~\ref{X8} and Table~\ref{talist} yield that $\Phgs(\bar H_\pu(\Phi'_8;\Q))$ is $(1+\Ll^2t^{-3})(1+\Ll^3t^{-5})\Ll^{-8}t^{16}\cdot (\Ll^{-2}t^8+\Ll^{-1}t^7)\inv$.

We compute the Borel--Moore homology of $\Psi:=\Phi'_8\cup\Phi'_9$ by exploiting the long exact sequence
\begin{equation}\label{psi}
\cdots\rightarrow\bar H_k(\Phi'_8;\Q)\rightarrow\bar H_k(\Psi;\Q)\rightarrow\bar H_k(\Phi'_9;\Q)\xrightarrow{\delta_k}\bar H_{k-1}(\Phi'_8;\Q)\rightarrow\cdots
\end{equation}

Both the Borel--Moore homology of $\Phi'_8$ and $\Phi'_9$ are tensor products of the Borel--Moore homology of $\PGL(3)$. The question is whether their structure as tensor products of $\bar H_\pu(\PGL(3);\Q)$ is respected by the maps in~\eqref{psi} or not.

We computed in \S~\ref{VGI} that the strata $1$--$7$ do not contribute to the Borel--Moore homology of $F'_{10}$. Recall that $F'_8\cup F'_9$ is a vector bundle of rank $1$ over $\Psi$. By comparing the geometry of $\mathcal D$ and its projectivization, we can conclude that 
\begin{equation}\label{F810}
\bar H_\pu(\bigcup_{i=8}^{10}F'_i;\Q)\cong \bar H_{\pu}(\Psi;\Q)\otimes\bar H_\pu(\C^*;\Q).
\end{equation}

Isomorphism \eqref{icLH}, together with the computation of the Borel--Moore homology of $\bigcup_{i=1}^7F'_i$ in \S~\ref{VGI}, yields that the Borel--Moore homology of $\bigcup_{i=8}^{10}F'_i$ is a tensor product of $\bar H_\pu(\GL(3);\Q)\cong\bar H_\pu(\C^*;\Q)\otimes\bar H_\pu(\PGL(3);\Q)$. 
In view of \eqref{F810}, this property implies that the Borel--Moore homology of $\Psi$ is a tensor product of $\bar H_\pu(\PGL(3);\Q)$. The only possibility for this is that the maps of the exact sequence~\eqref{psi} respect the structure of $\Phi'_8$ and $\Phi'_9$ as tensor products of $\bar H_\pu(\PGL(3);\Q)$. 
This is important, because it implies that all differentials $\delta_k$ in~\eqref{psi} are determined, once one knows the rank of 
$$\delta_{25}\co  \bar H_{25}(\Phi'_9;\Q)\cong\bigoplus^2\Inv\otimes\Q(2)\rightarrow\bar H_{24}(\Phi'_8;\Q)\cong\Inv\otimes\Q(2).$$ 

We claim that $\delta_{25}$ has rank one. This would yield part~\eqref{een} in the claim. Note that, in view of \eqref{F810}, equality \eqref{een} implies \eqref{twee}.

Define $B\subset X'_9$ as the locus of configurations $(p,q,\{r_i\})$ such that $$p,q\notin Sing\left(\bigcup_ir_i\right),\ \ \ pq\notin\{r_1,r_2,r_3,r_4\}.$$ Denote by $\mathcal B\rightarrow B$ the restriction of the bundle  $\Phi'_9\rightarrow X'_9$ to $B$. 
Next, consider the locus $A\subset X'_8$ of configurations $(p,q,\{e_1,e_2,e_3,e_4\})$ such that $p\in e_1e_3$, $q\in e_2e_4$, $\{p,q\}\cap(\{e_1,e_2,e_3,e_4\}\cup(e_1e_2\cap e_3e_4)=\emptyset$. Denote by $\mathcal A\rightarrow A$ the restriction of the bundle $\Phi'_8\rightarrow X'_8$ to $A$.
Note that for every element $a=(p,q,\{e_i\})$ of $A$, the configuration 
$C_a:=(p,q,\{e_1e_2,e_1e_3,e_2e_4,e_3e_4\})$ is an element of $B$. This means that the face of the 4-dimensional open simplex lying above $a$ is identified (in $\ba\Lambda'\ba$) with one of the external faces of the 5-dimensional simplex contained in $\mathcal B$ which lies above $C_a\in B$. Moreover, this 4-dimensional open simplex is the only face of $C_a$ which lies in $\mathcal A$. Recall that the Borel--Moore homology of the union of an open simplex and one of its open faces is trivial, because  of the characterization of Borel--Moore homology as the relative homology of the one-point compactification of a space modulo the added point. This implies that the Borel--Moore homology of $\mathcal A\cup\mathcal B$ is trivial. 

We have the following chains of inclusions:
$$
\begin{array}{r@{}c@{}l}
\mathcal A &\underset{\text{open}}{\hookrightarrow} &\Phi'_8\\
\text{\scriptsize closed } \cap & & \cap\text{\scriptsize{ closed}}\\ 
\mathcal A\cup\mathcal B&\underset{\text{open}}{\hookrightarrow} &\Psi\\
\text{\scriptsize open } \cup & & \cup\text{\scriptsize{} open }\\ 
\mathcal B&\underset{\text{open}}{\hookrightarrow} &\Phi'_9.\\
\end{array}
$$

In particular, if we consider the long exact sequence in Borel--Moore homology associated to the closed inclusion $\mathcal A\hookrightarrow\mathcal A\cup\mathcal B$, we have that the map $\bar H_{25}(\mathcal B;\Q)\rightarrow\bar H_{24}(\mathcal A;\Q)$ is an isomorphism. By the computation of the Borel--Moore homology of $X'_8$ in Lemmas~\ref{X8} and~\ref{utile} we have that $\bar H_{24}(\mathcal A;\Q)\cong\bar H_{24}(\Phi'_8;\Q)$ and $\bar H_{25}(\mathcal B;\Q)\subset\bar H_{25}(\Phi'_9;\Q)$ are both one-dimensional.
\end{proof}

Now that the contribution of strata $8$--$10$ to the Borel--Moore homology of $\mathcal D$ is known, we want to deduce their contribution to the rational cohomology of $\I_2$. Then the closed inclusion $\mathcal D\rightarrow\mathcal V$ (for the definition, see the beginning of the present section) induces a long exact sequence
$$\cdots\rightarrow \coh k{\mathcal V}\rightarrow\coh k{\I_2}\rightarrow \bar H_{33-k}(\mathcal D;\Q)(-k)\rightarrow\coh {k+1}{\mathcal V}\rightarrow\cdots,$$
which in the case of the part of Borel--Moore homology of $\mathcal D$ that comes from strata $8$--$10$, gives the following contribution to the $\s_2$-HG polynomial of the rational cohomology of $\I_2$:
\begin{equation}\label{secondblock1}
(1+\Ll t)(1+\Ll^2 t^3)(1+\Ll^3 t^5)\big(\Ll^6t^6\inv+\Ll^7t^8(\inv+\ant)+\Ll^8t^8\ant\big),
\end{equation}
hence the contribution of these strata to the HG polynomial of the cohomology of $\mathcal Q_2$ is
\begin{equation}\label{secondblock}
\Ll^6t^6\inv+\Ll^7t^8(\inv+\ant)+\Ll^8t^8\ant.
\end{equation}

\section{Vassiliev--Gorinov's method}\label{VGmethod}
In order to make the article as self-contained as possible, we include here an introduction to Vassiliev--Gorinov's method for computing the cohomology of complements of discriminants, following \cite{OTM4} and \cite{OTtesi}. This review of the method is by no means complete, and we encourage the interested reader to consult \cite{Vart}, \cite{Gorinov} and \cite{OTM4}.

Let $Z$ be a projective variety, $\mathcal F$ a vector bundle on $Z$ and $V$ the space of global sections of $\mathcal F$. Define the discriminant $\Sigma\subset V$ as the locus of sections with a vanishing locus which is either singular or not of the expected dimension. Assume that $\Sigma$ is a subvariety of $V$ of pure codimension $1$. 
Our aim is to  compute the rational cohomology of the complement of the discriminant, $X=V\setminus \Sigma$. This is equivalent to determining the Borel--Moore homology of the discriminant, because there is an isomorphism between the reduced cohomology of $X$ and Borel--Moore homology of $\Sigma$. If we denote by $M$ the dimension of $V$, this isomorphism can be formulated as
$$\tilde H^{\pu}(X;\Q)\cong \bar H_{2M-\pu-1}(\Sigma;\Q)(-M).$$

\begin{dfn}
A subset $S\subset Z$ is called a \emph{configuration} in $Z$ if it is compact and non-empty. The space of all configurations in $Z$ is denoted by $\Conf(Z)$.
\end{dfn}

\begin{prop}[\cite{Gorinov}]
The Fubini--Study metric on $Z$ induces in a natural way on $\Conf(Z)$ the structure of a compact complete metric space.
\end{prop}

To every element in $v\in V$, we can associate its singular locus $K_v\in\Conf(Z)\cup\{\emptyset\}$. We have that $K_0$ equals $Z$, and that $L(K):=\{v\in V: K\subset K_v\}$ is a linear space for all $K\in \Conf(Z)$. 

Vassiliev--Gorinov's method is based on the choice of a collection of families 
of configurations $X_1,\ldots,X_R\subset\Conf(Z)$, satisfying some axioms (\cite[3.2]{Gorinov}, \cite[List 2.1]{OTM4}). Intuitively, we have to start by classifying all possible singular loci of elements of $V$. Note that singular loci of the same type have a space $L(K)$ of the same dimension. We can put all singular configurations of the same type in a family. Then we order all families we get according to the inclusion of configurations. In this way we obtain a collection of families of configurations which may already satisfy Gorinov's axioms. If this is not the case, the problem can be solved by adding new families to the collection. Typically, the elements of these new families will be degenerations of configurations already considered. 
For instance, configurations with three points on the same projective line and a point outside it can degenerate into configurations with four points on the same line, even if there is no $v\in V$ which is only singular at four collinear points.  

Once the existence of a collection $X_1,\dots,X_R$ satisfying Gorinov's axioms is established, Vassi\-liev--Gorinov's method gives a recipe for constructing a space $\ba \mathcal X\ba$ and a map
$$|\epsilon|\co  \ba\X\ba\longrightarrow \Sigma,$$
called \emph{geometric realization}, which is a homotopy equivalence and induces an isomorphism on Borel--Moore homology.
The original construction by Vassiliev and Gorinov uses topological joins to construct $\ba \X\ba$. This construction was reformulated in \cite{OTM4} by using the language of cubical spaces. This ensures in particular that the map induced by $|\epsilon|$ on Borel--Moore homology respects mixed Hodge structures.

Vassiliev--Gorinov's method provides also a stratification $\{F_j\}_{j=1,\dots,N}$ on $\ba\X\ba$. Each $F_j$ is locally closed in $\ba\X\ba$, hence one gets a spectral sequence converging to $\bar H_\pu(\Sigma;\Q)\cong\bar H_\pu(\ba\X\ba;\Q)$, with $E^1_{p,q}\cong\bar H_{p+q}(F_p)$. To compute the Borel--Moore homology of $F_j$ for all $j=1,\dots,R$, it is helpful to use an auxiliary space $\ba\Lambda\ba$, whose construction depends only on the geometry of the families $X_1\dots,X_R$, and which is covered by locally closed subsets $\{\Phi_j\}_{j=1,\dots,N}$. 

\begin{prop}[\cite{Gorinov}]\label{ucci}
\renewcommand{\labelenumi}{\arabic{enumi}.}
\begin{enumerate}
\item\label{ucpri} For every $j=1,\dots,R$, the stratum $F_j$ is a complex vector bundle  over $\Phi_j$. The space $\Phi_j$ is in turn a fiber bundle over the configuration space $X_j$. 
\item\label{ucsec} If $X_j$ consists of configurations of $m$ points, the fiber of $\Phi_j$ over any $x\in X_j$ is an $(m-1)$-dimensional open simplex, which changes its orientation under the homotopy class of a loop in $X_j$ interchanging a pair of points in $x_j$.
\item\label{opeco}  If $X_R=\{Z\}$, $F_R$ is the open cone with vertex a point (corresponding to the configuration $Z$),  over $\ba\Lambda\ba\setminus\Phi_R$. 
\end{enumerate}\end{prop}

We recall here the topological definition of an open cone.
\begin{dfn}
Let $B$ be a topological space. Then a space is said to be an \emph{open cone} over $B$ with vertex a point if it is homeomorphic to the space $B\times[0,1)/R$, where the equivalence relation is $R=(B\times\{0\})^2$.
\end{dfn}

The fiber bundle $\Phi_j\rightarrow X_j$ of Proposition \ref{ucci} is in general non-orientable. As a consequence, we have to consider the homology of $X_j$ with coefficients not in $\Q$, but in some local system of rank one. Therefore we recall some constructions concerning Borel--Moore homology of configuration spaces with twisted coefficients.

\begin{dfn}\label{tuples}
Let $Z$ be a topological space. Then for every $k\geq 1$ we have the space of ordered configurations of $k$ points in $Z$,
$$\F kZ=Z^k\setminus\bigcup_{1\leq i<j\leq k}\{(z_1,\dots,z_k)\in Z^k: z_i=z_j\}.$$
There is a natural action of the symmetric group $\s_k$ on $F(k,Z)$. The quotient is called the space of unordered configurations of $k$ points in $Z$,
$$\B kZ=\F kZ/\s_k.$$
\end{dfn}

The \emph{sign representation} $\pi_1(\B kZ)\rightarrow \Aut(\Z)$ maps the paths in $\B kZ$ defining odd (respectively, even) permutations of $k$ points to multiplication by $-1$ (respectively, 1). The local system $\pm\Q$ over $\B kZ$ is the one locally isomorphic to $\Q$, but with monodromy representation equal to the sign representation of $\pi_1(\B kZ)$. We will often call $\bar H_\pu (\B kZ,\pm\Q)$ the \emph{Borel--Moore homology of $\B kZ$ with twisted coefficients}, or, simply, the \emph{twisted Borel--Moore homology of $\B kZ$}.

The following is Lemma~2 in \cite{Vart}.
\begin{lem}\begin{enumerate}\label{lem1}
\item
If $N\geq 1$, $k\geq2$, the twisted Borel--Moore homology of $\B k{\C^N}$ is trivial.
\item
If $N\geq1$, we have isomorphisms $$\bar H_\pu(\B k{\Pp N};\pm\Q)\cong\bar H_{\pu-k(k-1)}(\bb G(k-1,\Pp N);\Q)$$ for every $k\geq1$, where $\bb G(k-1,\Pp N)$ denotes the Grassmann variety of $(k-1)$-dimensional linear subspaces in $\Pp N$. In particular, $\bar H_\pu(\B k{\Pp N};\pm\Q)=0$ if $k\geq N+2$.
\end{enumerate}
\end{lem}

\def\cprime{$'$}

\end{document}